\documentclass[reqno]{amsart}

\usepackage{amsmath, amsthm, amscd, amsfonts, amssymb, graphicx, color}

\usepackage[bookmarksnumbered, colorlinks, plainpages]{hyperref}

\usepackage{lipsum}%% a garbage package you don't need except to create examples.
\usepackage{fancyhdr}

\usepackage{algorithm,algorithmic}

\newtheorem{theorem}{Theorem}[section]
\newtheorem{lemma}[theorem]{Lemma}
\newtheorem{proposition}[theorem]{Proposition}
\newtheorem{corollary}[theorem]{Corollary}
\theoremstyle{definition}

\theoremstyle{remark}
\newtheorem{remark}[theorem]{Remark}
\numberwithin{equation}{section}
\begin{document}

\title{A Refinement of the arithmetic-geometric mean inequality and some more }
\author{Mehdi Eghbali Amlashi}
%\address{Department of Pure Mathematics, Ferdowsi University
%of Mashhad, P. O. Box 1159, Mashhad 91775, Iran;}
%\curraddr{Department of ...}
%\email{id\_mail1@um.ac.ir}

\author{Mahmoud Hassani}

\address{Department of Mathematics, Mashhad Branch, Islamic Azad University, Mashhad, Iran;}

\email{Amlashi@mshdiau.ac.ir (Mehdi Eghbali Amlashi)}

\email{mhassanimath@gmail.com, hassani@mshdiau.ac.ir (Mahmoud Hassani)}
%\address{Department of Pure Mathematics, Ferdowsi University
%of Mashhad, P. O. Box 1159, Mashhad 91775, Iran;}
%\email{id\_mail2@um.ac.ir}

%\author{write your name}
%\address{Department of Pure Mathematics, Ferdowsi University
%of Mashhad, P. O. Box 1159, Mashhad 91775, Iran;}
%\email{id\_mail3@um.ac.ir}

%\thanks{thanks.}

%\date{}

%\subjclass[2000]{aa2bb; cc3dd.}

%\keywords{Graph Partition, Conjugate Gradient, Branch and bound algorithm
%(at least 3 and at most 5 items).}


\begin{abstract}
In this note, we present a refinement of the well-known AM-GM inequality. We use this improved inequalty to establish corresponding  inequalities on  Hilbert space. We also give some refinements of the Kantorovich inequality. 
\end{abstract}

\maketitle

\section{Introduction}

For any two positive real numbers a, b, their arithmetic and geometric means introduce as:
\[A(a, b) = \frac{a + b}{2}, \quad G(a, b) = \sqrt{ab}\] 
the classical AM-GM inequality says:
\[G(a, b) \leq A(a, b),\]
throught the years, some improvements have presented for this inequality. 
In \cite {7}, authors presented the following:
\begin{equation}\label{eq1.1}
\sqrt{ab} \leq \sqrt{ab} (1 + \frac{1}{8} (\ln a - \ln b)^2) \leq \frac{a + b}{2}.
\end{equation}
In the next proposition, we prove \ref{eq1.1}, by an elementary and simple method.
\begin{proposition}
For any two positive real numbers a and b:
\[\sqrt{ab} (1 + \frac{1}{8} (\ln (\frac{a}{b}))^2) \leq \frac{a + b}{2}.\]
\end{proposition}
\begin{proof}
We consider the function, $f: (0 ,\infty) \rightarrow \mathbb{R} $, defined by \\ $$f(x) = \sqrt{x} (1 +\frac{1}{8}(\ln x)^2).$$Then  $y = \frac{x +1}{2}$ is the tangant line of it's graph at $x = 1$.\\
Furthermore, by taking the second derivative of the function $f$, we get \\ $$ {f''}(x) = \frac{-1}{32x \sqrt{x}}(\ln x)^2.$$\\
 So, $f$ is concave on interval $(0, \infty)$. This implies: 
\begin{equation}\label{eq1.5}
\sqrt{x}(1 +\frac{1}{8}(\ln x)^2) \leq \frac{x +1}{2}.
\end{equation}
Now, put $x = \frac{a}{b}$. This completes the proof.
\end{proof}
The Heinz mean for two positive real numbers a, b is:
\[H_{\nu}(a, b) = \frac{1}{2}(a^{1 - \nu}b^{\nu} + a^{\nu} b^{1 - \nu}), 0 \leq \nu \leq 1.\]
This is closely connected with the arithmetic and geometric means and it interpolates the two mentioned. So if $0 \leq \nu \leq 1$ we have
\[\sqrt{x} \leq \frac{1}{2}(x^{1 - \nu} + x^{\nu}) \leq \frac{x +1}{2},    for  x > 0.\]
By integrating w.r.t $\nu$ on $[0, 1]$,
\begin{equation}\label{eq1.2}
\sqrt{x} \leq \frac{x - 1}{\ln x} \leq \frac{x + 1}{2}, \quad x > 0
\end{equation}
(Recall that $\mathcal{L} (a, b) = \frac{a - b}{\ln a  - \ln b}$ is the logarithmic mean of positive real numbers a, b.)\\

\begin{theorem}
For any real numbers $ a, b > 0 $,
\begin{equation}\label{eq2.1}
\sqrt{ab} \leq \sqrt{ab}(1 + \frac{1}{2}(\frac{a - b}{a + b})^2) \leq \frac{a + b}{2}.
\end{equation}
\end{theorem}
\begin{proof}
By \ref{eq1.2}
\[\frac{\vert x - 1 \vert}{\vert \ln x \vert} = \frac{x - 1}{\ln x} \leq \frac{x + 1}{2}, \quad x > 0.\]
So, we have
\[\frac{\vert x - 1 \vert}{x + 1} \leq \frac{1}{2} \vert \ln x \vert .\]
And,
\[\left( \frac{x - 1}{x +1} \right)^2 \leq \frac{1}{4} (\ln x)^2 .\]
So by \ref{eq1.5} we get
\begin{equation}\label{eq2.7}
 \sqrt{x}(1 + \dfrac{1}{2} \left( \frac{x - 1}{x +1} \right)^2) \leq \sqrt{x}(1 + \frac{1}{8} (\ln x)^2) \leq \frac{x + 1}{2}, \quad x > 0. 
\end{equation}
Now \eqref{eq2.1} is proved by \eqref{eq2.7}, and putting $x = \frac{a}{b}$.
\end{proof}

\begin{remark}
	Note that the function $y = \frac{x - 1}{x + 1}$ is increasing on $(1, +\infty)$. So, for $1 < m < x $, we have
	\[\frac{m - 1}{m + 1} \leq \frac{x - 1}{x + 1} .\]
	Now,by \ref{eq2.7} we can write the following inequality:
	\begin{equation}\label{eq3.32}
		\sqrt{x} \left( 1 + \frac{1}{2}(\frac{m - 1}{m + 1})^2 \right) \leq \frac{x + 1}{2}, \quad 1 < m < x,
	\end{equation}
and if $0<ma<b$ , taking $x=\dfrac{b}{a}$ in  inequality (1.6), we have 
\begin{equation}\label{eq2.22}
	\sqrt{ab} \leq \sqrt{ab}(1 + \frac{1}{2}(\frac{m - 1}{m + 1})^2) \leq \frac{a + b}{2}.
\end{equation}
\end{remark}
\begin{theorem}
	For any real numbers $ x > 0 $,
	\begin{equation}\label{eq2.8}
		\dfrac{x + 1}{2} \leq \sqrt{x}(1 + \frac{1}{8}(\sqrt{x}-\dfrac{1}{\sqrt{x}})^{2}) .
	\end{equation}
\end{theorem}
\begin{proof}
	Consider  the function, $f: (0 ,\infty) \rightarrow \mathbb{R} $, defined by 
	\[ f(x) = \sqrt{x}(1 + \frac{1}{8}(\sqrt{x}-\dfrac{1}{\sqrt{x}})^{2}). \]
	By taking the first and second derivatives of the function $f$, we get 
	\[ f'(x)= \dfrac{3}{16}\sqrt{x} - \dfrac{1}{16x\sqrt{x}} + \dfrac{3}{x\sqrt{x}},\]
	and,
	\[ f''(x) = \dfrac{3}{32\sqrt{x}} + \dfrac{3}{32x^{2}\sqrt{x}}- \dfrac{3}{16x\sqrt{x}} \] 
	\[= \dfrac{3}{32x^{2}\sqrt {x}}(x^{2}-2x +1)\geq 0. \]
	So, $f$ is convex on interval $(0, \infty)$. But, $y=\dfrac{x+1}{2}$ is the tangent line of $y=f(x)$ at $x=1$. Therefore $\dfrac{x+1}{2} \leq f(x)$, and the proof is complete.	
\end{proof}

\section{Operator inequalities for positive linear maps}
In this section, with the help of Theorem 1.2., we deduce refinement of the arithmetic- geometric mean inequality for Hilbert space operators. Let $H$ be a complex Hilbert space. We represent the set of all bounded operators on $H$ by $B(H)$.
If  $A\in B(\textit{H})$ satisfies $A^{*}=A$, then $A$ is called a self-adjoint operator. If a self-adjoint operator $A$ satisfies $\langle x, Ax \rangle \geq 0 $ for any $ x  \in \textit{H}$ , then $A$ is called a positive operator . For two self-adjoint operators $A$ and $B$ , $A \geq B$ means $A-B \geq 0$. The notation $A > 0$ means $A$ is an invertible positive operator.
A linear map $\Phi $ on $B(H)$ is positive if $\Phi(A) \geq 0$  whenever $A \geq 0$. It is said to be unital if $\Phi(I) = I$.   
 Let A  and B be two positive operators in $B(\mathcal{H})$. Arithmetic mean of A and B has familiar form:
\[\frac{A + B}{2}.\]
For $ A,B > 0 $, the geometric mean $A  \sharp B$ is defined by \\
\[A \sharp B = A^{\frac{1}{2}}(A^{-\frac{1}{2}} B A^{- \frac{1}{2}})^{\frac{1}{2}} A^{\frac{1}{2}}.\]
Indeed, $A \sharp B$ is the unique positive solution to the Riccati equation:
\[X A^{-1} X = B.\]
One motivation of such a notion is of course the AM-GM inequality:
\[A \sharp B \leq \frac{A + B}{2}.\]
\begin{theorem}\label{thm2.2}
	For any positive invertible operators A and B:
	\[A \sharp B \leq (A \sharp B) \left( I + \frac{1}{2} ((A+B)^{-1}(A-B))^2 \right) \leq \frac{A+B}{2}.\]
\end{theorem}

\begin{proof}
By \ref{eq2.7}
\[x^{\frac{1}{2}} \leq x^{\frac{1}{2}} (1 + \frac{1}{2}((x+1)^{-1}(x-1))^2) \leq \frac{x+1}{2}.\]
 By applying a standard functional calculus for the positive operator $T = A^{-\frac{1}{2}} B A^{-\frac{1}{2}}$, we infer from the previous inequality,
\[T^{\frac{1}{2}} \leq T^{\frac{1}{2}} + \frac{1}{2}T^{\frac{1}{2}} \left( (T + I)^{-1}(T -I)\right)^2 \leq \frac{T + I}{2}.\]
Multiplying $A^{\frac{1}{2}}$ to the above inequality from left-hand side and right-hand side, we have 
\[A \sharp B \leq A \sharp B + \frac{1}{2}M \leq \frac{A+B}{2},\]
where

\[ M = A^{\frac{1}{2}} T^{\frac{1}{2}} \left( (T + I)^{-1}(T -I)\right)^2 A^{\frac{1}{2}} \]
\[ = A^{\frac{1}{2}} T^{\frac{1}{2}} A^{\frac{1}{2}} \left( A^{- \frac{1}{2}} (T + I)^{-1} A^{- \frac{1}{2}}\right) \left( A^{\frac{1}{2}} (T - I) A^{\frac{1}{2}} \right)
\left( A^{- \frac{1}{2}} (T + I)^{-1} A^{- \frac{1}{2}}\right) \left( A^{\frac{1}{2}} (T - I) A^{\frac{1}{2}} \right) \]
\[ = (A \sharp B) \left( (A + B)^{-1}(A - B)\right)^2. \]

\end{proof}
\begin{theorem}\label{thm2.2}
	For any positive invertible operators $A$ and $B$ we have 
	\[ \dfrac{A+B}{2} \leq \dfrac{1}{8}(A\sharp B)( A^{-1}B+ B^{-1}A + 6I).\]
\end{theorem}

\begin{proof}
	By \ref{eq2.8}
	\[ \dfrac{x + 1}{2} \leq x^{\frac{1}{2}}( \frac{1}{8}x + \dfrac{1}{8}x^{-1} + \dfrac{3}{4}) .\]
	We use functional calculus with the positive operator $T = A^{-\frac{1}{2}} B A^{-\frac{1}{2}}$, we get
	\[\dfrac{T+1}{2} \leq T^{\frac{1}{2}}( \frac{1}{8}T + \dfrac{1}{8}T^{-1} + \dfrac{3}{4}I). \]
	Multiplying $A^{\frac{1}{2}}$ to the above inequality from left-hand side and right-hand side, we have 
	\[ \dfrac{A+B}{2} \leq A^{\frac{1}{2}}T^{\frac{1}{2}}A^{\frac{1}{2}}( \dfrac{1}{8}A^{\frac{-1}{2}}TA^{\frac{1}{2}}+\dfrac{1}{8}A^{\frac{-1}{2}}T^{-1}A^{\frac{1}{2}} + \dfrac{3}{4}I ).\]
	Hence 
		\[ \dfrac{A+B}{2} \leq \dfrac{1}{8}(A\sharp B)( A^{-1}B+ B^{-1}A + 6I).\].
\end{proof}

\section{Refinements of kantorovich inequality}

The  Kantorovich inequality can be stated as follows \cite {5}.\\
 For an operator $A$ such that 
$0 < mI  \leq A \leq M I$, where $m, M$ are real numbers and for any unit vector $x$ :

\[\langle Ax, x \rangle \langle A^{-1}x, x \rangle \leq \frac{(m + M)^2}{4 m M }. \]
This inequality was interesting for many mathematicians and some refinements of  Kantorovich inequality were presented.  Before we give the main results in this section,  let us present the following Lemma that will be useful later. 
\begin{lemma}
	Let $A, B \in B(\mathcal{H})$ and $1 <m<M$ such that  
	$ 0 < m B \leq A \leq M B$ 
	. Then 
	\[ 	\left( 1 + \frac{1}{2}(\frac{m - 1}{m + 1})^2 \right) A \sharp B \leq \frac{A + B}{2}.
	\]
\end{lemma}
\begin{proof}
	
	Since $mI \leq B^{-\frac{1}{2}}AB^{-\frac{1}{2}}\leq MI$ by using functional calculus, we infer from inequality(1.6), 
	\[(B^{-\frac{1}{2}}AB^{-\frac{1}{2}})^{\frac{1}{2}} \left( 1 + \frac{1}{2}(\frac{m - 1}{m + 1})^2 \right) \leq \frac{B^{-\frac{1}{2}}AB^{-\frac{1}{2}} + I}{2}. \]
	Multiplying both sides by $B^{\frac{1}{2}}$ on the left and right, we deduce the desired result.

\end{proof}

In the following we obtain an improvement of Kantorovich inequality.
\begin{theorem}
Let $A, B \in B(\mathcal{H})$ such that  
$ 0 < mI \leq m' B \leq A \leq M I$ 
. Then for every unit vector $x \in H$,
\[\langle Ax, x \rangle \langle B x, x \rangle \leq \frac{(m + M)^2}{4 m M }(1 + \frac{1}{2} (\frac{m' - 1}{m' + 1})^2)^{-2} \langle (A \sharp B)x, x \rangle^{2}. \]
\end{theorem}
\begin{proof}
	
	According to the hypothesis  we get the order relation,
	\[\frac{m m'}{M} I \leq B^{-\frac{1}{2}} A B^{-\frac{1}{2}} \leq \frac{M m'}{m} I,\]
	and by $T =( B^{-\frac{1}{2}} A B^{-\frac{1}{2}})^{\frac{1}{2}}$:
	\[ 0 \leq \left (T - \sqrt{\frac{m m'}{M}} I \right) \left( \sqrt{\frac{m' M}{m}} I - T \right),\]
	which yeilds:
	\[m' I + T^{2} \leq  T( \sqrt{\dfrac{Mm'}{m}} + \sqrt{\dfrac{mm'}{M}}). \]
	By multiplying $B^{\frac{1}{2}}$ to left and right of this inequality:
	\begin{equation}\label{eq3.33}
	A +m' B \leq ( \sqrt{\dfrac{Mm'}{m}} + \sqrt{\dfrac{mm'}{M}}) A \sharp B.
		\end{equation}
	So,
	\[ \langle Ax, x \rangle + m' \langle Bx, x \rangle \leq  ( \sqrt{\dfrac{Mm'}{m}} + \sqrt{\dfrac{mm'}{M}})\langle (A \sharp B)x, x \rangle.\]
	But, by \ref{eq3.32} 
	\[2 \sqrt{m' \langle Ax, x \rangle \langle Bx, x \rangle} (1 + \frac{1}{2} (\frac{m' - 1}{m' + 1})^2) \leq  \langle Ax, x \rangle + m' \langle Bx, x \rangle.\]
	
	This completes the proof.
\end{proof}
If we choose $B =A ^{-1}$ we get from Theorem 3.2 the following Corollary. 
\begin{corollary}
		Let $A \in B(\mathcal{H})$ be a positive operator such that, 
	\[I < mI < m' A^{-1} \leq A \leq M I,\]
	then
	\begin{equation}\label{eq3.17}
	\langle Ax, x \rangle \langle A^{-1}x, x \rangle \leq \frac{(m + M)^2}{4 m M }(1 + \frac{1}{2} (\frac{m' - 1}{m' + 1})^2)^{-2}.
	\end{equation}
\end{corollary}
 Inequality (3.2) provides a refinement of Kantorovich inequality.\\
In 1996, using the operator geometric mean, Nakamoto and Nakamura \cite {6} proved that
\begin{equation}\label{eq3.5}
	 \Phi(A) \sharp \Phi(A^{-1}) \leq \dfrac{M+m}{2\sqrt{Mm}},
\end{equation}
whenever $0<mI<A<MI$ and $\Phi $ is a unital positive linear map on $B(H)$.\\
Our second main result in this section, which is related to Inequality (3.3) can be stated as follows:
 \begin{theorem}
 Let $\Phi $ be a unital positive linear map on $ B(H)$ and	let $A, B \in B(\mathcal{H})$ such that  
 	$ 0 < mI \leq m' B \leq m'^{2} B \leq A \leq M I$ 
 	. Then for every unit vector $x \in H$,
 	\[\Phi(A) \sharp \Phi(B) \leq \frac{m + M}{2\sqrt{Mmm'} }(1 + \frac{1}{2} (\frac{m' - 1}{m' + 1})^2)^{-1}  \Phi (A \sharp B).\]
 \end{theorem}
\begin{proof}
	According to the hypothesis we have
	\[   0 < mI \leq m'(\Phi(m'B)) \leq \Phi (A) \leq M , \] 
	and based on Lemma 1.3, we get the order relation,
	
	\[	\left( 1 + \frac{1}{2}(\frac{m' - 1}{m' + 1})^2 \right) \phi(A) \sharp \phi(m'B) \leq \frac{\phi(A) +\phi(m' B)}{2}. \]

		But, by \ref{eq3.33}
	\[ \Phi(A)  + m' \Phi(B)  \leq  ( \sqrt{\dfrac{Mm'}{m}} + \sqrt{\dfrac{mm'}{M}})  \Phi(A \sharp B).\]
	So
	\[ 2 \sqrt{ m'}\left( 1 + \frac{1}{2}(\frac{m' - 1}{m' + 1})^2 \right) \phi(A) \sharp \phi(B) \leq ( \sqrt{\dfrac{Mm'}{m}} + \sqrt{\dfrac{mm'}{M}})  \Phi(A \sharp B). \]
This completes the proof.
	\end{proof}
If we choose $B =A ^{-1}$ we get from Theorem 3.4 that
 \begin{corollary}
	Let $\Phi $ be a unital positive linear map on $ B(H)$ and	let $A\in B(\mathcal{H})$ such that  
	$ 0 < mI \leq m'A^{-1} \leq m'^{2}A^{-1}\leq A \leq M I$ 
	. Then for every unit vector $x \in H$,
	
	\begin{equation}\label{eq3.21}
	\Phi(A) \sharp \Phi(A^{-1}) \leq \frac{m + M}{2\sqrt{Mmm'} }(1 + \frac{1}{2} (\frac{m' - 1}{m' + 1})^2)^{-1}. 
		\end{equation}
\end{corollary}
 Inequality (3.4) provides a refinement of Inequality (3.3).

\end{document}